


\magnification\magstep1
\baselineskip=14pt
\vsize=24.0truecm
\font\bigbf=cmbx12
\font\csc=cmcsc10
\def\hatt{\widehat}
\def\tilda{\widetilde}
\def\half{\hbox{$1\over2$}}
\def\eps{\varepsilon}
\def\dd{{\rm d}}
\def\E{{\rm E}}

\def\midd{\,|\,}
\def\data{{\rm data}}

\def\cstok#1{\leavevmode\thinspace\hbox{\vrule\vtop{\vbox{\hrule\kern1pt
        \hbox{\vphantom{\tt/}\thinspace{\tt#1}\thinspace}}
        \kern1pt\hrule}\vrule}\thinspace} 
\def\firkant{\cstok{\phantom{$\cdot$}}}

\centerline{\bigbf Bayesian and Empirical Bayesian Bootstrapping} 

\medskip
\centerline{\bf Nils Lid Hjort}
\medskip
\centerline{\sl University of Oslo and Norwegian Computing Centre} 

\smallskip
\centerline{\sl -- May 1991 --}

\smallskip

{{\smallskip\narrower 
\baselineskip=12pt\noindent
{\csc Abstract.} 
Let $X_1,\ldots,X_n$ be a random sample from an unknown probability
distribution $P$ on the sample space ${\cal X}$, and let 
$\theta=\theta(P)$ be a parameter of interest. The present paper
proposes a nonparametric `Bayesian bootstrap' method of obtaining Bayes 
estimates and Bayesian confidence limits for $\theta$.
It uses a simple simulation technique to numerically 
approximate the exact posterior distribution of $\theta$ 
using a (non-degenerate) Dirichlet process prior for $P$. 
Asymptotic arguments are given which justify the use of the 
Bayesian bootstrap for any smooth functional $\theta(P)$. 
When the prior is fixed and the sample size grows five approaches 
become first-order equivalent: the exact Bayesian, the Bayesian
bootstrap, Rubin's degenerate-prior bootstrap, Efron's bootstrap,
and the classical one using delta methods. 
The Bayesian bootstrap method is also extended to the
semiparametric regression case. 
A separate section treats similar ideas for 
censored data and for more general hazard rate models,
where a connection is made to a `weird bootstrap' proposed by Gill. 
Finally empirical Bayesian versions of the procedure are discussed,
where suitable parameters of the Dirichlet process prior 
are inferred from data. 

Our results lend Bayesian support to the classic Efron bootstrap.
It is the Bayesian bootstrap under a noninformative reference prior; 
it is a limit of natural approximations 
to good Bayes solutions; 
it is an approximation 
to a natural empirical Bayesian strategy; 
and the formally incorrect reading of a bootstrap histogram as
a posterior distribution for the parameter isn't so incorrect after all.

\smallskip }}

\bigskip

{\bf 1. Introduction and summary.}
Let $X_1,\ldots,X_n$ be independent and identically distributed 
(i.i.d.)~according to an unknown distribution $P$. 
For convenience take the sample
space to be ${\cal X}={\cal R}$, the real line, so that $P$ can be 
identified with its distribution function (c.d.f.)~$F$.
Most of the methods and results in this report have natural 
extensions to ${\cal R}^k$ and indeed to any complete,
separable metric space $\cal X$.

Let $\theta=\theta(F)$ be a parameter functional of interest,
like the mean, or median, or the standard deviation,
or $\int|x-F^{-1}(\half)|\,\dd F(x)$. 
We shall be concerned with nonparametric Bayesian estimates of 
and confidence statements about $\theta$,
and need to start out with a prior distribution
on the space of all c.d.f.'s. A natural class from which to choose is 
provided by Ferguson's (1973, 1974) Dirichlet processes;
this class is rich, each member has large support, 
and at least posterior expectations 
(Bayes estimates under quadratic loss) 
can be calculated for a fair list of cases. Thus let 
$$F\sim {\rm Dir}(aF_0),\eqno(1.1)$$
i.e.~$F$ is a Dirichlet process with parameter $aF_0$. 
Here $F_0(.)=\E_B\,F(.)$ is the prior guess c.d.f.~whereas $a>0$ 
has interpretation as prior sample size, see Ferguson (op.~cit.).
Subscript $B$ means that the operation in question 
is relative to the Bayesian framework. 

The observed sample $x_1,\ldots,x_n$ gives rise to and can
be identified with the empirical 
c.d.f.~$F_n(t)={1\over n}\sum_{i=1}^nI\{x_i\le t\}$. 
The posterior distribution of $F$ is
$${\cal L}_B\{F\midd\data\}
	={\cal L}_B\{F\midd F_n\}={\rm Dir}(aF_0+nF_n).\eqno(1.2)$$
Thus the distribution function
$$G_n(t)={\rm Pr}_B\{\theta(F)\le t\midd \data\}
	={\rm Pr}_B\{\theta(F)\le t\midd F_n\} \eqno(1.3)$$
is in principle known to the statistician. 
In addition to the Bayesian point estimate
$\theta_B=\E_B\,\{\theta(F)\midd \data\}=\int t\,\dd G_n(t)$
we wish to calculate Bayesian confidence limits 
$\theta_L$ and $\theta_U$ from the data, obeying
${\rm Pr}_B\{\theta_L\le\theta(F)\le\theta_U\midd\data\}\doteq1-2\alpha$, 
say. Thus
$$\theta_L=G_n^{-1}(\alpha),\qquad \theta_U=G_n^{-1}(1-\alpha)\eqno(1.4)$$
are the natural choices, where 
$G_n^{-1}(p)=\inf\{t\colon G_n(t)\ge p\}$.

The fact that $G_n$ above is only very rarely explicitly available, 
however, necessitates devising computational approximations. 
There are a couple of rather laborious ways to simulate 
variables from a close approximation to $G_n$.
That is, a sequence $Y_1,Y_2,\ldots$ being i.i.d.~with a distribution
very close to $G_n$ can be generated, thus enabling one 
to obtain a close approximation to $G_n$ and to the sought-after
$\theta_B$, $\theta_L$, $\theta_U$; see 2A and 9B below. 
It turns out that the following rather simpler 
alternative simulation strategy gives  
a good approximation to the posterior distribution $G_n$: 
Generate a `Bayesian bootstrap (BB) sample' 
$X_1^*,\ldots,X_{n+a}^*$ of size $n+a$ from the mixture distribution 
$$F_{n,B}(t)=\E_B\,\{F(t)\midd \data\}
	={a\over a+n}F_0(t)+{n\over a+n}F_n(t),\eqno(1.5)$$
the natural Bayes estimate of the underlying c.d.f.~$F$, and compute
a `BB parameter value' 
$$\theta_{BB}^*=\theta(F_{BB}^*)
	=\theta\Bigl({1\over n+a}\sum_{i=1}^{n+a}I\{X_i^*\le\,.\,\}\Bigr)
	=\theta(X^*_1,\ldots,X^*_{n+a}) \eqno(1.6)$$
on the basis of the empirical c.d.f.~$F_{BB}^*$ of these values.
The proposed approximation to $G_n$ is 
$$\hatt G_n(t)={\rm Pr}_*\{\theta_{BB}^*\le t\},\eqno(1.7)$$
where subscript `$_*$' is used to indicate operations relative to
the (data-conditional) BB framework. 
In practice $\hatt G_n(.)$ is evaluated via simulations as 
$$\hatt G_n(t)
	\doteq\hatt G_{n,{\rm boot}}(t)
	={1\over {\rm boot}}\sum_{b=1}^{\rm boot}
		I\{\theta^{*b}_{BB}\le t\}$$
for a large number ${\rm boot}$ of independent $\theta_{BB}^{*b}$
of the type described. 
This idea, inserting $\hatt G_n$ for $G_n$, 
leads to using 
$$\hatt\theta_B=\int t\,\dd\hatt G_n(t)
  \doteq {1\over{\rm boot}}\sum_{b=1}^{\rm boot}\theta^{*b}_{BB}, 
	\eqno(1.8)$$
the BB-based Bayes estimate of $\theta(F)$, 
and to the BB percentile interval
$$\hatt\theta_L=\hatt G_n^{-1}(\alpha)\doteq\hatt G_{n,\rm boot}(\alpha)
	\le \theta(F)\le
  \hatt\theta_U=\hatt G_n^{-1}(1-\alpha)
	\doteq\hatt G_{n,\rm boot}(1-\alpha). \eqno(1.9)$$
Other Bayesian posterior calculations can be carried out with
the same relative ease, like computing Bayes estimates with
non-quadratic loss functions. 

The motivation for the BB method lies in the fact that the two conditional
distributions ${\cal L}_B\{F\midd \data\}$ and
${\cal L}_*\{F^*_{BB}\midd \data\}$ are reasonably close. 
This is explained in Section~2.
When the prior sample size parameter $a$ goes to zero 
(corresponding in a certain sense to the case of a 
`noninformative nonparametric prior' for $F$) 
the BB becomes Efron's classic bootstrap.
In particular the BB percentile interval becomes 
Efron's (uncorrected) percentile interval. 
One may therefore think of the BB as an `informative extension'
of the usual bootstrap method, capable of incorporating
prior information on $F$. 
This also lends some Bayesian credit to Efron's bootstrap,
and shows that the incorrect interpretation of the 
traditional bootstrap distribution as a posterior distribution for
the parameter isn't so incorrect after all. 
Note that the BB smooths also outside the data points, 
unlike the classic bootstrap.  
In the $a$ close to zero case the BB is also an approximation 
to Rubin's (1981) `degenerate prior Bayesian bootstrap',
as indicated in Section 2. 
A large-sample justification for the BB is given in Section 3. 
Under frequentist circumstances it is shown that 
five different approaches tend to agree, asymptotically;
the exact Bayesian, 
the bootstrap Bayesian, 
the Rubin bootstrap, 
the ordinary nonparametric large-sample method, 
and the classic bootstrap.   

Section 4 gives two Bayesian bias correction methods for the BB. 
In Section 5 the BB method is shown at work for 
a couple of parameters. 
Some suggestions on how to select 
parameters in the prior Dirichlet process is briefly 
discussed in Section 6, thus opening for empirical Bayes versions
of the bootstrap. 
In particular the Rubin method, for which the Efron method 
is the BB approximation, can be seen as a natural empirical Bayes
strategy.  
Section 7 presents the BB for semiparametric regression,
where the residual distribution is given a Dirichlet prior.
This leads in particular to an interesting frequentist 
bootstrap scheme suggestion. 
In Section 8 we deviate a bit from the main story and 
report on a brief investigation into bootstrapping schemes 
for censored data and hazard rate models. 
Finally several supplementing remarks are made in Section 9. 


\bigskip
{\bf 2. The Bayesian bootstrap.} 
This section motivates the Bayesian bootstrap method (1.5)--(1.9)
and explains why it can be expected to work. 
Then connections to Efron's (1979, 1982) traditional bootstrap and 
to Rubin's (1981) degenerate-prior Bayesian bootstrap are commented on.

\smallskip
{\sl 2A. Approximating the posterior distribution.} 
Before discussing our Bayesian bootstrap method further, 
let us mention that the posterior distribution $G_n$ of (1.3)
can be evaluated exactly for a few parameter functionals. 
Section 5 provides calculations for 
$\theta=F\{A\}$, $A$ a set of interest, 
and $\theta=F^{-1}(p)$, the $p$-quantile. 
For other parameters it may be possible to carry out 
almost-exact simulation of ${\cal L}_B\{\theta\midd\data\}$, 
as hinted at before (1.5). 

For such an example, let $\theta=\int x\,\dd F(x)$ be the mean of $F$.
The exact distribution of $\theta$ given data can be obtained,
but the resulting expressions are complicated and make exact simulation
difficult. See Hannum, Hollander, and Langberg (1981), Yamato (1984),
and Cifarelli and Regazzini (1990).
However, the posterior distribution can be approximated 
with that of $\theta'=\sum_{i=1}^nx_iF\{x_i\}+\sum_{j=1}^my_jF\{A_j\}$,
say, where $A_1,\ldots,A_m$ is a fine partition of 
${\cal R}-\{x_1,\ldots,x_n\}$, and $y_j\in A_j$. This $\theta'$ can then
be simulated, since $(F\{x_1\},\ldots,F\{x_n\},F(A_1),\ldots,F(A_m))$
has a (finite-dimensional) Dirichlet distribution. 
Hjort (1976) showed that $\beta_m\rightarrow_d\beta$ in $\cal X$
implies ${\rm Dir}(\beta_m)\rightarrow_d{\rm Dir}(\beta)$
in the space of probability measures on $\cal X$, 
w.r.t.~various metrics, and 
$\int x\,\dd F_m(x)\rightarrow_d\int x\,\dd F(x)$ under a mild 
extra condition on $\{\beta_m\}$. 
This result justifies ${\cal L}(\theta)\approx{\cal L}(\theta')$ above,
and can be used to approximate $G_n$ also in more general cases, 
using a simpler variable that involves only 
finite-dimensional Dirichlet distributions. 
Another almost-exact simulation strategy is described in Section 9B. 

This example illustrates that (1.3) in general will be difficult
to obtain via exact or almost-exact simulation from $G_n$. 
The Bayesian bootstrap method described in (1.5)--(1.9)
is clearly much simpler. Note that 
$X^*_i$ is from $F_0$ with probability $a/(a+n)$ and 
is equal to $x_j$ with probability $1/(a+n)$, for $j=1,\ldots,n$. 
The description in (1.6)--(1.7) assumed $a$ to be an integer. 
If $a=m+\beta$, say, with $m$ an integer and $0<\beta<1$, generate 
$n+m+1$ $X^*_i$'s from $F_{n,B}$ instead, and use 
$F_{BB}^*(t)=\bigl[\sum_{i=1}^{n+m}I\{X_i^*\le t\}
		+\beta I\{X^*_{n+m+1}\le t\}\bigr]/(n+m+\beta)$. 

To explain why the BB method can be expected to work,
consider the two data-conditional distributions 
${\cal L}_B\{F\midd \data\}$ and ${\cal L}_*\{F_{BB}^*\midd \data\}$. 
Judicious calculations give 
$$\eqalign{
\E_B\,\{F(t)\midd\data\}&=F_{n,B}(t), \cr
\E_*\,\{F_{BB}^*(t)\midd\data\}&=F_{n,B}(t), \cr
{\rm cov}_B\bigl[\{F(s),F(t)\}\midd\data\bigr]
	&={1\over n+a+1}F_{n,B}(s)\{1-F_{n,B}(t)\}, \cr
{\rm cov}_*\bigl[\{F_{BB}^*(s),F_{BB}^*(t)\}\midd\data\bigr]
	&={1\over n+a}F_{n,B}(s)\{1-F_{n,B}(t)\}, \cr} \eqno(2.1)$$
for all $s\le t$.
Accordingly, for well-behaved functionals $\theta=\theta(F)$ we would expect
$${\cal L}_B\{\theta(F)\midd\data\}\approx 
	{\cal L}_*\{\theta(F_{BB}^*)\midd\data\}, 
	\quad {\rm that\ is\ }G_n\doteq \hatt G_n. \eqno(2.2)$$
As a point of further comparison it may be considered a bit annoying that 
the skewness of $F\midd\data$ is about twice that of
$F_{BB}^*\midd\data$, but they are both small for moderate to large~$n$:
$$\eqalign{
\E_B\,\{F(t)-F_{n,B}(t)\}^3\midd \data
	&={2F_{n,B}(t)\{1-F_{n,B}(t)\}\{1-2F_{n,B}(t)\}
				\over (n+a+1)(n+a+2)}, \cr
\E_*\,\{F_{BB}^*(t)-F_{n,B}(t)\}^3\midd \data
    	&={F_{n,B}(t)\{1-F_{n,B}(t)\}\{1-2F_{n,B}(t)\}
				\over (n+a)^2}. \cr} \eqno(2.3)$$
One might therefore expect the (uncorrected) 
BB and the exact Bayes methods to
be first order but not second order equivalent; see Section 3. 

We could have made the second order moments agree completely
and not only approximatively in (2.1) by drawing $n+a+1$ BB-data,
instead of $n+a$, to form $F_{BB}^*$. 
The difference is tiny and disappears for moderate to 
large samples. We have chosen BB-sample size $n+a$
to better reflect the Bayesian balancing of prior information
and data and to better highlight the generalisation from 
the usual Efron bootstrap. 

\smallskip
{\sl 2B. Connections to other bootstraps.}
Consider the non-informative case $a$ close to zero
(or, rather, $a/n$ close to zero). Then the BB procedure 
advocates taking bootstrap samples of size $n$ from the usual $F_n$,
and basing analysis on simulating $\theta^*=\theta(X^*_1,\ldots,X^*_n)$. 
But this is the familiar nonparametric Efron bootstrap! 
In particular the BB percentile interval becomes Efron's 
(uncorrected) percentile interval. 
Thus the BB method is a proper Bayesian generalisation of the
classic bootstrap. 
And since the BB really works, by Section~3, 
this also lends Bayesian credit to the classic bootstrap;
it is the `vague prior' version of a natural nonparametric Bayesian
strategy. The incorrect interpretation of the bootstrap 
distribution (say in the form of a histogram of 1000 
bootstrap values) as a posterior distribution for the parameter
isn't that incorrect after all; 
it is an approximation to the true posterior distribution 
if the starting point is a Dirichlet with a small $a$. 

There are better confidence interval
methods than the percentile method for the classic bootstrap, 
but the more sophisticated versions, incorporating bias and 
acceleration corrections, are still first-order large-sample
equivalent to the simple one. Corrections to the BB percentile
interval appear in Section 4 below.  
It should also be remarked that 
the classic bootstrap has several other uses than
the making of confidence intervals, like estimating 
variances of complicated estimators. 
The BB scheme is general enough to handle such problems too,
but would in general need an inside bootstrap loop as well. 

Rubin (1981) and Efron (1982, Ch.~10) discuss 
a simple Bayesian bootstrap different from
the one proposed here. 
The Rubin bootstrap, although somewhat differently presented 
in his paper, can be seen to be the limiting Bayes method obtained 
by using $F\sim{\rm Dir}(aF_0)$ as prior and 
then letting $a\rightarrow0$, i.e.~using 
${\cal L}_B\{\data\}={\rm Dir}(nF_n)$, see (1.2). 
(Actually, Rubin and Efron consider only finite sample spaces,
but the extension to the present generality is not difficult 
using the available theory of Dirichlet processes.)
In this limiting case $F\midd\data$ is concentrated on the observed 
data values, $F=\sum_{i=1}^nd_i\delta(x_i)$, with
weights $(d_1,\ldots,d_n)$ following a Dirichlet $(1,\ldots,1)$
distribution (uniform on the simplex of nonnegative weights
summing to one). 
In particular values of $\theta(F)$ 
can be simulated from the exact $G_n$ of (1.3).
The $d_i$'s may be simulated as $e_i/(e_1+\cdots+e_n)$,
where the $e_i$'s are unit exponential. 
If $\theta(F)=\int x\,\dd F(x)$ is the mean, 
for example, then a large number of realisations 
of $\theta(F)=\sum_{i=1}^nd_ix_i=\sum_{i=1}^ne_ix_i/\sum_{i=1}^ne_i$
can be generated, the distribution of these values will approximate $G_n$,
enabling one to get good numerical approximations to $\theta_B$ and
to the interval (1.4). 
Rubin (1981) notes that this approach, 
though different in interpretation, agrees well,
operationally and inferentially, with the ordinary bootstrap procedure. 

The Rubin bootstrap does not come out of letting $a\rightarrow0$
in the BB method proposed here. Results of Section 3 show that 
these are large-sample equivalent to the first order,
in particular the Bayesian using a Dirichlet prior with a small $a$ 
(Rubin) may view the Efron bootstrap (which is our BB with small $a$)
as a numerical simulation device giving approximately the same results. 
Rubin's method smooths the weights, but rigidly sticks to the 
observed data poins (as does the ordinary bootstrap), 
whereas the more generally applicable BB method proposed here smooths also 
outside the data points, using $F_{n,B}$.
One might call this paper's BB the informative Bayesian bootstrap
and Rubin's BB the degenerate-prior bootstrap. 
(And with due fairness Rubin didn't advocate its general use,
but concentrated on connections to and comparisons with 
Efron's method.)  
The results and remarks above suggest 
that the present informative BB 
comes much closer to being a proper Bayesian generalisation of 
Efron's bootstrap, both in operation and in spirit.

In a recent paper Newton and Raftery (1991) have developed a
Bayesian-inspired weighted likelihood bootstrap. 
In its nonparametric form it generalises Rubin's method in 
a way different from our BB. Their method does not smooth outside
the data points, whereas our does, in presence of a prior guess $F_0$. 
See further discussion in their Section 8.
There are finally indirect connections 
to some of the bootstrapping schemes we discuss
for hazard rate models in Section 8 below. 

\bigskip
{\bf 3. Large-sample justification: Five statisticians agree.}
In this section it is proved that 
the two conditional distributions ${\cal L}_B\{\theta(F)\midd\data\}$
and ${\cal L}_*\{\theta(F)\midd \data\}$ are asymptotically
equivalent to the first order. We also show that five 
different approaches tend to give the same inference for large samples;
the classical using delta methods, 
the classic bootstrap,
the accurate Bayesian using Dirichlet priors, 
Rubin's non-informative prior bootstrap, 
and the BB. 
Then some supplementing remarks are made.

\smallskip
{\sl 3A. Finite sample space.} 
Assume first, and mostly for illustrational purposes, 
that the sample space is finite, say ${\cal X}=\{1,\ldots,L\}$. Let 
$$f_{\rm true}(l)={\rm Pr}_F\{X_i=l\}, \quad 
  f_n(l)={1\over n}\sum_{i=1}^nI\{x_i=l\}, \quad 
  {\rm and\ }f_{n,B}(l)={af_0(l)+nf_n(l)\over a+n}.$$ 
Efron (1982, Ch.~5.6) observed that 
$${\cal L}\{\sqrt{n}(f_n-f_{\rm true})\}
	\rightarrow N_L\{0,\Sigma(f_{\rm true})\}, \eqno(3.1)$$
$${\cal L}_*\{\sqrt{n}(f^*_n-f_n)\midd\data\}
	\approx N_L\{0,\Sigma(f_n)\}
	\rightarrow N_L\{0,\Sigma(f_{\rm true})\}{\rm \ a.s.}, \eqno(3.2)$$
where $f^*_n(l)=(1/n)\sum_{i=1}^nI\{\tilda x_i=l\}$ 
stems from the ordinary bootstrap, and where 
$\Sigma(f)$ has elements 
$f(l)\delta_{l,m}-f(l)f(m)$. 
Efron discussed why (3.1)--(3.2) 
may be taken as an asymptotic justification
for a class of inferential procedures based on the bootstrap. 
Note that the (3.1)--(3.2) results rely 
only on asymptotic theory for the multinomial distribution, 
and that the `almost surely' statement refers to the set $\Omega_0$
of probability 1 under which each $f_n(l)$ goes to $f_{\rm true}(l)$.  

These can now be accompanied by results for the exact and the BB 
approximated posterior distributions 
${\cal L}_B\{f\midd\data\}$, ${\cal L}_*\{f_{BB}^*\midd\data\}$.
The framework is the frequentist one,
where the $X_i$'s are truly i.i.d.~from $f_{\rm true}$. 
One can prove 
$${\cal L}_B\bigl\{(n+a+1)^{1/2}(f-f_{n,B})\midd\data\bigr\}
	\approx N_L\{0,\Sigma(f_{n,B})\}
	 \rightarrow N_L\{0,\Sigma(f_{\rm true})\}{\rm \ a.s.}, \eqno(3.3)$$
$${\cal L}_*\bigl\{(n+a)^{1/2}(f_{BB}^*-f_{n,B})\midd\data\bigr\}
	\approx N_L\{0,\Sigma(f_{n,B})\}
	 \rightarrow N_L\{0,\Sigma(f_{\rm true})\}{\rm \ a.s.} \eqno(3.4)$$
The first follows from asymptotic properties of the Dirichlet distribution,
while the second is essentially the multidimensional central limit theorem.
Note that exactly the same a.s.~set $\Omega_0$ is at work in (3.1)--(3.4). 
The parameter $a$ is supposed to be fixed in (3.3)--(3.4), 
so that $f_{n,B}\rightarrow f_{\rm true}$ on $\Omega_0$,  
but arguments underlying the indicated approximations 
show that the two distributions are approximately equal 
even if $a$ goes to infinity with~$n$. 
[Certain minimax procedures correspond to using $a$ proportional to 
$\sqrt{n}$, for example; see Hjort (1976).]
Efron's discussion of the consequences 
of (3.1) and (3.2) (1979, p.~23; 1982, Ch.~5.6)
can now be applied to (3.3) and (3.4) as well, 
and provides the asymptotic justification 
for the BB procedure for the case of a finite sample space. 

It is interesting to note that if only $a/\sqrt{n}\rightarrow0$
as $n$ grows, then $\sqrt{n}\{f(l)-f_n(l)\}-\sqrt{n}\{f(l)-f_{n,B}(l)\}$
goes to zero, which implies 
$${\cal L}_B\{\sqrt{n}(f-f_n)\midd\data\}
	\rightarrow N_L\{0,\Sigma(f_{\rm true})\}
	{\rm \ a.s.}, \eqno(3.5)$$
$${\cal L}_*\{\sqrt{n}(f_{BB}^*-f_n)\midd\data\}
	\rightarrow N_L\{0,\Sigma(f_{\rm true})\} {\rm \ a.s.} \eqno(3.6)$$ 
Accordingly, looking back at (3.1)--(3.2), 
four different approaches will lead to the same 
inferential statements, up to first order asymptotics:
the classical based on $f_n$;
the ordinary Efron bootstrap;
the proper posterior Bayes;
and the BB. This holds for each fixed $a$, also for $a\rightarrow0$,
which means that Rubin's degenerate-prior bootstrap (see Section 2) 
also is large-sample equivalent to the other four. 

\smallskip
{\sl 3B. The real line.} 
Now consider extension of the preceding results and conclusions 
to ${\cal X}={\cal R}$. The degree to which (3.1) and (3.2) and its 
consequences have analogues for ${\cal X}={\cal R}$ was 
investigated in Bickel and Freedman (1981) and Singh (1981), 
and later on in the form of extensions and refinements 
by other authors. The canonical parallel to (3.1) is
$${\cal L}\bigl[\sqrt{n}\{F_n(.)-F(.)\}\bigr]
	\rightarrow W_0(F(.))
	\quad {\rm in\ }D[-\infty,\infty], \eqno(3.7)$$
where $W_0$ is a Brownian bridge
and convergence takes place in the space $D[-\infty,\infty]$
of all right continuous functions $y(.)$ on the line with left hand limits and 
obeying $y(-\infty)=y(\infty)=0$, 
see for example Billingsley (1968). 
Bickel and Freedman (1981) proved the bootstrap companion
$${\cal L}_*\bigl[\sqrt{n}\{F_n^*(.)-F_n(.)\}\midd\data\bigr]
	\rightarrow W_0(F(.))
	\quad {\rm in\ }D[-\infty,\infty] {\rm \ a.s.}, \eqno(3.8)$$
and concluded that the bootstrap works for well-behaved functionals
$\theta=\theta(F)$. 

These results can be parallelled in the present Bayesian posterior context. 
Again, we look at limiting properties in an ordinary framework 
in which $F_n$ according to the Glivenko--Cantelli theorem 
converges uniformly to $F=F_{\rm true}$ 
on a set $\Omega_0$ of probability one. 

\smallskip
{\csc Theorem.} 
{\sl Let $a$ vary with $n$ in such a way that 
$F_{n,B}=(aF_0+nF_n)/(a+n)$ goes to some $F_\infty$ on $\Omega_0$;
$F_\infty$ is just $F_{\rm true}$ if only $a/n$ goes to
zero. Then
$${\cal L}_B\bigl[(n+a+1)^{1/2}\{F(.)-F_{n,B}(.)\}\midd\data\bigr]
	\rightarrow W_0(F_\infty(.)), \eqno(3.9)$$
$${\cal L}_*\bigl[(n+a)^{1/2}\{F_{BB}^*-F_{n,B}(.)\}\midd\data\bigr]
	\rightarrow W_0(F_\infty(.)), \eqno(3.10)$$
along every sequence in $\Omega_0$.} 

{\csc Proof:}
The second statement is within reach of the (triangular version of)
the classical Donsker invariance theorem for i.i.d.~random variables. 
The first statement involves proving finite-dimensional convergence 
and tightness. Finite-dimensional convergence follows upon studying 
asymptotic properties of (finite-dimensional) Dirichlet distributions. 
To prove tightness (with probability 1), 
calculate first $\E\,(k+1)^2(U-\alpha)^2(V-\beta)^2$ where 
$(U,V,W)$ is Dirichlet $(k\alpha,k\beta,k\gamma)$ and $\alpha+\beta+\gamma=1$. 
The resulting expression can be bounded by $3\alpha\beta$, 
regardless of $k$. Hence 
$$\eqalign{
(n+a+1)^2\,\E_B\,\bigl[F(s,t]&-F_{n,B}(s,t]\}^2
	\{F(t,u]-F_{n,B}(t,u]\}^2\midd\data\bigr] \cr
	&\le 3F_{n,B}(s,t]F_{n,B}(t,u] \cr}$$
for $s\le t\le u$. Taking limsup gives the bound 
$3F_\infty(s,t]F_\infty(t,u]$ on the right hand side,
for sequences in $\Omega_0$. 
This implies tightness by the proof of 
Billingsley's (1968) Theorem 15.6 
(but not quite by the theorem itself). \firkant

\smallskip
Thus the conditional distributions $\theta(F)\midd\data$
and $\theta(F_{BB}^*)\midd\data$ will be close 
to each other for well-behaved functionals, justifying the BB method.
Particular examples can be worked through, 
as in Bickel and Freedman (1981). 
Their tentative description of well-behavedness (p.~1209) 
can also be subscribed to here. Sufficient conditions for
$${\cal L}_B\bigl[(n+a+1)^{1/2}\{\theta(F)-\theta(F_{n,B})\}\midd\data\bigr]
	\rightarrow N\{0,\sigma^2(F_\infty)\} {\rm\ a.s.},$$
$${\cal L}_*\big[(n+a)^{1/2}\{\theta(F_{BB}^*)
    	-\theta(F_{n,B})\}\midd\data\bigr]
	\rightarrow N\{0,\sigma^2(F_\infty)\} {\rm\ a.s.}$$
to hold, for appropriate variance $\sigma^2(F_\infty)$, 
can be written down using von Mises or influence function methods.
See for example Boos and Serfling (1980) and Parr (1985),
who use Frech\'et differentiability, 
or Shao (1989) who uses Lipschitz differentiability, 
or Gill (1989) with Hadamard or compact differentiability. 
The limit results obtainable using such machinery imply 
$$\hatt G_n^{-1}(p)\doteq G_n^{-1}(p)
	\doteq\theta(F_{n,B})+z_p\,\sigma(F_\infty)/\sqrt{n+a},$$
where $z_p$ is the $p$-quantile of the standard normal. 

If $a$ is fixed, or only $a/\sqrt{n}\rightarrow0$, then 
$${\cal L}_B\bigl[\sqrt{n}\{F(.)-F_n(.)\}\midd \data\bigr]
	\rightarrow W_0(F_{\rm true}(.)){\rm \ a.s.}, $$
$${\cal L}_*\bigl[\sqrt{n}\{F_{BB}^*(.)-F_n(.)\}\midd \data\bigr]
	\rightarrow W_0(F_{\rm true}(.)){\rm \ a.s.} $$
A conclusion concerning the approximate agreement among the 
five statisticians referred to after (3.5)--(3.6) is 
therefore reached also for ${\cal X}={\cal R}$ 
(and for more general spaces). 
Each of them reaches confidence intervals that are first-order
equivalent to 
$$\theta(F_n)-z_{1-\alpha}\,\sigma(F_n)/\sqrt{n}
	\le \theta(F) \le 
  \theta(F_n)+z_{1-\alpha}\,\sigma(F_n)/\sqrt{n}, \eqno(3.11)$$
albeit from different perspectives and with partly 
different interpretations. This holds for each ${\rm Dir}(aF_0)$ 
prior, and, regarding the classic bootstrap, 
holds for both the simple percentile interval and for the somewhat
better reflected bootstrap interval 
$[2\theta(F_n)-G^*_n{}^{-1}(1-\alpha),2\theta(F_n)-G^*_n{}^{-1}(\alpha)]$,
where $G^*_n$ is the bootstrap distribution. 

It is perhaps surprising that a simple method like the BB,
constructed merely to make the mean function and covariance function
of the exact and approximate distributions of $F(.)$ agree, 
can work well for the vast majority of parameter functionals. 
As indicated in (3.9)--(3.10) 
this is at least partly the work and the magic of the 
functional central limit theorem. 
This also points to the possibility of 
using `small-sample asymptotics' machinery to arrive at other
approximations to the posterior distribution $G_n$,
for example Edgeworth--Cram\'er expansions combined with 
Taylor expansions. Such an approach would be functional-dependent,
however; a primary virtue of the BB 
is that it is both simple and versatile. 
A similar remark applies to the classic bootstrap, of course. 


The results in this section are taken from the technical report 
Hjort (1985). Results resembling (3.9) and (3.10) 
have also been found by Lo (1987), who also worked with 
rates of convergence.  

\bigskip
{\bf 4. A Bayesian bias correction to the BB percentile interval.}
The ordinary frequentist bootstrap percentile intervals can be corrected for
bias and acceleration, see Efron (1987). 
The BB percentile interval (1.9) cannot be corrected in the same way,
cf.~Hjort (1985, Section 4). 
There is another possibility of detecting and repairing a bias,
however. For each in a respectable catalogue of examples there
is a known transformation $h$, perhaps the identity, such that the 
posterior expected value of $h(\theta(F))$ is explicitly 
calculable by some published formula, 
i.e.~$\nu_0=\E_B\,\{h(\theta(F))\midd \data\}$ is known. 
The BB procedure uses 
$$\hatt H_n(t)={\rm Pr}_*\{h(\theta(F^*_{BB}))\le t\midd \data\}
	=\hatt G_n(h^{-1}(t))$$
to estimate $H_n$, the c.d.f.~of $h(\theta(F))$ given data, 
and approximates $\nu_0$ with 
$$\hatt\nu_0=\int t\,\dd\hatt H_n(t)\doteq {1\over {\rm boot}}
	\sum_{i=1}^{\rm boot}h(\theta(F^{*b}_{BB}))=\nu_0+\eps, \eqno(4.1)$$
say. Accordingly, if $\eps\not=0$, then $\hatt H_n$ 
is not a perfect estimate of $H_n$. The repaired estimate
$\hatt H_\eps(t)=\hatt H(t+\eps)$ gets the mean right, however. Hence
$$\hatt H_\eps^{-1}(\alpha)
	=\hatt H^{-1}(\alpha)-\eps\le h(\theta(F))
	\le \hatt H^{-1}(1-\alpha)-\eps=\hatt H_\eps(1-\alpha)$$
would be a natural corrected confidence interval for $h(\theta(F))$.
Transforming back we obtain 
$$h^{-1}\bigl[h(\hatt G_n^{-1}(\alpha))-\eps\bigr]
	\le \theta(F) \le
  h^{-1}\bigl[h(\hatt G_n^{-1}(1-\alpha))-\eps\bigr] \eqno(4.2)$$
as the {\sl bias-corrected BB percentile interval} for $\theta(F)$.
Of course this interval is just (1.9) if $\eps=0$. 
We emphasise that the bias correction is not concerned with
frequentist coverage probabilities, but is a simple way 
of repairing the BB estimate $\hatt G_n$ of $G_n$ so as to 
get the mean of $h(\theta(F))$ straight.  

As an example, suppose an interval is needed for $\sigma(F)$,
the standard deviation. One may prove, using methods 
of Ferguson (1973) and Hjort (1976), that 
$$\eqalign{
\E_B\,\{\sigma^2(F)\midd \data\}
&={n+a\over n+a+1}\Bigl[{a\over n+a}\sigma^2(F_0)
	+{n\over n+a}\sigma^2(F_n) \cr
&\qquad\qquad 	+{a\over n+a}{n\over n+a}\{\theta(F_n)
		-\theta(F_0)\}^2\Bigr]. \cr}\eqno(4.3)$$
The bias corrected confidence interval for $\sigma(F)$ is therefore
$$\{\hatt G_n^{-1}(\alpha)^2-\eps\}^{1/2}\le \sigma(F) \le
  		\{\hatt G_n^{-1}(1-\alpha)^2-\eps\}^{1/2},$$
where $\eps$ is the difference between the average value of the
observed $(\sigma^*_{BB})^2$ and $\E_B\,\{\sigma^2(F)\midd \allowbreak\data\}$. 

One can also write down a slightly more general 
{\sl bias and variance corrected} BB percentile interval which also 
takes into account the value of 
$\tau_0^2={\rm Var}\{h(\theta(F))\midd \data\}$ if it is available. 
Assume that, in addition to (4.1), 
$$\hatt\tau_0^2=\int (t-\hatt\nu_0)^2\,\dd\hatt H_n(t)
	\doteq{1\over {\rm boot}}\sum_{b=1}^{\rm boot}
	\{h(\theta^{*b}_{BB}-\hatt\nu_0\}^2=\tau_0^2(1+\delta)^2.$$ 
A perhaps better estimate of $G_n(h^{-1}(t))$ is then
$\hatt H_{n,\eps,\delta}(t)=\hatt H_n((1+\delta)t+\eps-\nu_0\delta)$,
since it manages to get both the mean and the variance right. 
Using $\hatt H_{n,\eps,\delta}^{-1}(p)
	=\{\hatt H_n^{-1}(p)+\nu_0\delta-\eps\}/(1+\delta)$ one ends up with
$$h^{-1}\Bigl[{h(\hatt G_n^{-1}(\alpha))+\nu_0\delta-\eps
		\over 1+\delta}\Bigr]
	\le\theta(F)\le 
  h^{-1}\Bigl[{h(\hatt G_n^{-1}(1-\alpha))
		+\nu_0\delta-\eps\over 1+\delta}\Bigr]. \eqno(4.4)$$
For an example, consider the mean parameter 
$\theta(F)=\int x\,\dd F(x)$, for which the posterior expectation is 
$\nu_0={a\over a+n}\theta(F_0)+{n\over a+n}\bar X_n$
and the posterior variance is 
$$\tau_0^2={1\over n+a+1}\Bigl[{a\over n+a}\sigma^2(F_0)
	+{n\over n+a}\sigma^2(F_n)
	+{a\over n+a}{n\over n+a}\{\theta(F_n)
		-\theta(F_0)\}^2\Bigr]. $$
The last formula is proved using methods of Ferguson (1973) and 
Hjort (1976) again, cf.~(4.3). 
The bias and variance corrected interval for $\theta(F)$ is 
$${\hatt G_n^{-1}(\alpha)+\nu_0\delta-\eps
		\over 1+\delta}
	\le\theta(F)\le 
  {\hatt G_n^{-1}(1-\alpha)+\nu_0\delta-\eps
		\over 1+\delta}. $$
One can similarly handle parameters of the type
$g\bigl(\int f(x)\,\dd F(x)\bigr)$. 

It should also be possible to construct a Bayesian 
skewness correction, cf.~(2.3), but this is not pursued here. 

\bigskip
{\bf 5. Some exact calculations.}
This section looks briefly into the nature of the BB approximation
method in two cases where exact calculations are possible.

\smallskip
{\sl 5A. A probability.}
If $\theta(F)=F(A)$ for some set $A$ of interest, then
$${\cal L}_B\{\theta(F)\midd \data\}
	={\rm Beta}\bigl\{aF_0(A)+\#(x_i\in A),
		a(1-F_0(A))+\#(x_i\notin A)\bigr\}.$$
Thus (1.3) and (1.4) can be obtained from tables of the incomplete
Beta function. In this case the BB method amounts to approximating 
the Beta distribution $G_n$ with that of $Y/(n+a)$, where 
$Y$ is binomial $[n+a,\{aF_0(A)+\#(x_i\in A)\}/(a+n)]$. 

If $U$ is ${\rm Beta}\{mp,m(1-p)\}$ and 
$V$ is ${\rm Bin}\{m,p\}/m$, then $\E\,U=\E\,V=p$, and 
$${\rm Var}\,U={p(1-p)\over m+1},
	\quad 
	{\rm Var}\,V={p(1-p)\over m}.$$
They differ in skewness and kurtosis, but not to any dramatic extent;
for example, 
$${\rm skew}\,U=2{(m+1)^{1/2}\over m+2}{1-2p\over \{p(1-p)\}^{1/2}},
	\quad
  {\rm skew}\,V={1\over m^{1/2}}{1-2p\over \{p(1-p)\}^{1/2}}.$$
Brief investigations have shown the distributions of $U$ and $V$,
and therefore confidence intervals based on either the 
exact or BB approximated distributions, to be remarkably similar, 
even for moderate $m$. This holds 
provided $p$ is not too close to zero or one,
provided $\alpha$ is not too close to zero, and finally provided 
the discrete distribution of $V$ is interpolated. 
Rather than using $\hatt G_m(t)={\rm Pr}[{\rm Bin}\{m,p\}/m\le t]$,
which jumps at the points $j/m$, use 
$\tilda G_m(j/m)=\half{\rm Pr}[{\rm Bin}\{m,p\}/m\le j/m]
+\half{\rm Pr}[{\rm Bin}\{m,p\}/m\le (j-1)/m]$, and interpolate linearly 
in between. Similar modifications to $\hatt G_n$ of (1.7)
should also be used  
in other cases where it increases in sharp jumps. 

\smallskip
{\sl 5B. The median.}
The $p$-quantile functional is another example where it is possible
to calculate the posterior distribution explicitly, but the 
resulting expressions are complex, and the BB would be much easier to
carry out in practice. For simplicity only the 
median $\theta(F)=F^{-1}(\half)=\inf\{t\colon F(t)\ge \half\}$ 
is considered here. 

Assume for concreteness that 
the data points are distinct, with $x_1<\ldots <x_n$.
We shall find $G_n(t)={\rm Pr}\{\theta(F)\le t\midd \data\}$.
For data point $x_j$ one has
$$\eqalign{
G_n\{x_j\}&={\rm Pr}\{F(-\infty,x_j)<\half,F(-\infty,x_j]\ge\half\} \cr
	&={\rm Pr}\{U<\half,U+V\ge\half\}
	 ={\rm Pr}\{U<\half,W<\half\}, \cr}$$
in which $(U,V,W)$ is Dirichlet with parameters 
$\alpha=aF_0(x_j-)+j-1$,
$\beta=aF_0\{x_j\}+1$,
and $\gamma=aF_0(x_j,\infty)+n-j$. Taking the prior guess c.d.f.~to be
continuous we find
$$\eqalign{
G_n\{x_j\}&={\Gamma(\alpha+\beta+\gamma)\over 
	\Gamma(\alpha)\Gamma(\beta)\Gamma(\gamma)}
	{1\over \alpha}\Bigl({1\over2}\Bigr)^\alpha	
	{1\over \gamma}\Bigl({1\over2}\Bigr)^\gamma \cr
&={\Gamma(a+n)\over \Gamma(aF_0(x_j)+j)\Gamma(a\{1-F_0(x_j)\}+n-j+1)}
	\Bigl({1\over2}\Bigr)^{a+n-1}. \cr}$$
Next consider $G_n[t,t+dt]$ for some $t$ outside the data points,
and let for further convenience $F_0$ be the integral of a prior guess density
$f_0$. Following the reasoning above one may show that 
$G_n$ has density at $t\in(x_j,x_{j+1})$ given by
$$g_n(t)={\Gamma(a+n)
	\over \Gamma(aF_0(t)+j)\Gamma(a\{1-F_0(t)\}+n-j)}
	\,af_0(t)\,J\bigl[aF_0(t)+j,a\{1-F_0(t)\}+n-j\bigr],$$
where
$$J[\alpha,\gamma]=\int_0^{\half}\int_0^{\half}
	u^{\alpha-1}w^{\gamma-1}(1-u-w)^{-1}\,\dd u\,\dd w. $$
It is in principle possible to compute for example the posterior
expectation and the upper and lower 5 percent points for $G_n$ based on
this.

Now consider the BB method in this situation.
Let for simplicity $n+a=2m+1$ be odd. 
The BB approximates the complicated $G_n$ using 
$X^*_i$'s from $F_{n,B}$, as follows:
$$\eqalign{
\hatt G_n(t)&={\rm Pr}_*\bigl[\theta^*_{BB}
	={\rm median}\{X^*_1,\ldots,X^*_{n+a}\}\le t\midd \data\bigr] \cr
&={\rm Pr}\bigl[{\rm Bin}\{2m+1,F_{n,B}(t)\}\ge m+1\bigr]. \cr}$$
Expressions for $\hatt G_n\{x_j\}$ and for the density $\hatt g_n$ 
that the distribution has between data points can be worked out based on this,
and they can be compared with $G_n$ and $g_n$ obtained above.
Such a study is not pursued here. Note that the endpoints of the
BB confidence interval (1.9) can be found using binomial tables. 
Note finally that in the non-informative case, where $a\rightarrow0$,
both $G_n$ and $\hatt G_n$ are supported on the data points. 

\bigskip
{\bf 6. Empirical Bayesian bootstrapping.}
The ideal Bayesian is able to specify $a$ and $F_0$ from the 
infamous but abstract `prior considerations'. 
Results of Section 3 show that the 
importance of these parameters diminishes and disappears 
with growing $n$, but they do matter for small and moderate $n$. 
This section briefly discusses some empirical methods. 

\smallskip
{\sl 6A. Choosing $a$ and parameters in $F_0$.} 
In some situations previous data may be available 
that are either of the same type as the $X_i$'s or at least
of a similar type. In the best case one has $m$ previous measurements
$X_i^0$ that come from the same $F$ as the new $X_i$'s. 
Then one may use $a=m$ (indeed the `prior sample size')
and $F_0$ equal to a smoothed empirical distribution or 
some fitted normal, say. 

In other cases one might have a
specified candidate $F_0$ from previous similar data, 
but without knowing for certain that the new data are from the
same distribution. Then the problem is to choose $a$, either
from informal `strength of belief' considerations, or from 
the new data. One wants to use a small $a$ if data disagree with 
the old $F_0$ and a larger one if they seem to fit. 
This can be done in a formal way by looking at moment properties
of the empirical distribution $F_n$. We have 
$\E\,(F_n-F_0)^2\midd F=(F-F_0)^2+F(1-F)/n$, so that 
$$\E\,(F_n-F_0)^2=\E\,(F-F_0)^2+{1\over n}\E\,F(1-F)
	=\Bigl({1\over a+1}+{1\over n}{a\over a+1}\Bigr)F_0(1-F_0),$$
and this can be used to fit a suitable $a$, for example via
$$\E\int(F_n-F_0)^2\,dW={1\over n}\Bigl(1+{n-1\over a+1}\Bigr)
		\int F_0(1-F_0)\,\dd W,$$
which holds for each weight measure $W$. Choosing $W=F_0$ gives 
$1/6$ for the last integral (in the continuous case). 

There are other estimation methods for $a$ with a fixed $F_0$ and
that to a larger extent uses properties of the Dirichlet process.
The maximum likelihood estimator can be derived, see Hjort (1976).
This and some other estimators depend however on the ties configurations
in the data in a somewhat strange way, and the sufficient statistic
is $D_n$, the number of distinct data points. 
This stems from some of the more esoteric and less 
satisfying mathematical properties of samples from a Dirichlet, 
and equating moments of 
natural quantities like above seems much more reasonable. 

In still other cases there might be parameters in the prior guess
$F_0$ to specify, say $F_0=N\{\mu_0,\sigma_0^2\}$. Then
one is helped by $\E\,\bar X_n=\E\int x\,\dd F(x)=\mu_0$ and 
$$\E\,{1\over n-1}\sum_{i=1}^n(X_i-\bar X_n)^2
	=\E\,\sigma^2(F)={a\over a+1}\sigma_0^2. $$
If $F_0$ is nonsymmetric one might fit parameters using also
$$\eqalign{
\E\,{n\over (n-1)(n-2)}\sum_{i=1}^n(X_i-\bar X_n)^3
	&=\E\int\{x-\mu(F)\}^3\,\dd F(x) \cr
	&={a\over a+1}{a\over a+2}\int\{x-\mu(F_0)\}^3\,\dd F_0(x). \cr}$$
All these moment methods should be used with care and sense. 
The moment formulae here have been 
proved using methods in Ferguson (1973) and Hjort (1976). 

We mention finally that a frequentist inspired 
double bootstrap method 
for fitting a good weight ${a\over a+n}$ in 
the mixture ${a\over a+n}F_0+{n\over a+n}F_n$ was suggested 
in Hjort (1988). 

\smallskip
{\sl 6B. Cross validation and Rubin--Efron as empirical Bayes solutions.}
Observe that the kind of schemes described above can be used on the basis of 
only the given data set, by dividing it into a small training set
and the remaining test set, or by some more elaborate
cross validation strategy. A simple version of this is as follows:
Pick $a$ data points to constitute the training set, 
from which the nonparametric guess on $F$ is $F_a$, the empirical
c.d.f.~for these. Since the remaining $n-a$ data points come from the
same $F$ the considerations above suggest using 
a ${\rm Dir}(aF_a)$ as prior for $F$. 
But then the posterior becomes Dirichlet with
$aF_a+(n-a)F_{n-a}=nF_n=\sum_{i=1}^n\delta(x_i)$. 
This is accordingly an empirical Bayes 
argument for using Rubin's method, and a fortiori for using 
its natural BB approximation, which is the classic Efron bootstrap.  

\bigskip
{\bf 7. BB in semiparametric regression.}
This section briefly discusses the extension of some of the previous
methods and results to the semiparametric regression case. 
The model is 
$$Y_i=\sum_{j=1}^px_{i,j}\beta_j+\sigma\eps_i
	=x_i'\beta+\sigma\eps_i, \quad i=1,\ldots,n, \eqno(7.1)$$
where the standardised residuals $\eps_i=(Y_i-x_i'\beta)/\sigma$
are i.i.d.~from $F$. The Bayesian version must have a prior distribution
for the unknown parameters $\beta$, $\sigma$, $F$. 
We stipulate that $(\beta,\sigma)$ comes from some prior density 
$\pi(\beta,\sigma)$ and that $F$, independently, comes from
the Dirichlet process with parameter $a\Phi$, where $\Phi(.)$
is the standard normal. When $a$ is large then the 
distribution of $F$ becomes concentrated in $\Phi$, 
which gives us the familiar textbook normal regression model.
This is accordingly a Bayesian generalisation with 
built-in uncertainty about the residual distribution. 
We are interested in Bayesian inference about parameters 
$\theta=\theta(\beta,\sigma,F)$, like regression deciles
$x'\beta+\sigma F^{-1}(j/10)$, probabilities 
${\rm Pr}\{Y(x)\le y\}=F((y-x'\beta)/\sigma)$, 
expected distance $\E\,|Y(x)-x'\beta|$, {\it\&}cetera. 
A BB strategy is arrived at below which makes it possible
to get an approximation to the full posterior distribution of 
such parameters. 

More general results of Hjort (1986, 1987) 
imply that the posterior density for $(\beta,\sigma)$ becomes 
$$\pi_n(\beta,\sigma\midd \data)={\rm const.}\,\pi(\beta,\sigma)
  \prod_{i=1}^n \bigl\{\sigma^{-1}
	\phi((y_i-x_i'\beta)/\sigma)\bigr\}, \eqno(7.2)$$
provided the $y_i$'s are distinct, i.e., the posterior
distribution for these parameters are as if $F$ had been known to be $\Phi$. 
And $F$ has a distribution being a mixture of Dirichlet processes, since
$$F\midd \{\beta,\sigma,\data\}\sim
	{\rm Dir}\bigl\{a\Phi+\sum_{i=1}^n
		\delta((y_i-x_i'\beta)/\sigma)\bigr\}. \eqno(7.3)$$
This makes it easy to write down $\E_B\,\{F(t)\midd \beta,\sigma,\data\}$ 
and then integrating out $(\beta,\sigma)\sim\pi_n(.)$ to reach
the posterior expectation of $F(t)$. 
Suppose for simplicity that $\sigma$ is known and that $\beta$ is given
a flat prior on ${\cal R}^p$, which leads to $\pi_n$ being 
quite simply $N\{\hatt\beta,{1\over n}\sigma^2M^{-1}\}$, where 
$M={1\over n}\sum_{i=1}^nx_ix_i'$ and 
$\hatt\beta=M^{-1}{1\over n}\sum_{i=1}^nx_iY_i$ 
is the familiar least squares estimator 
(now seen also as the Bayes solution under the flat prior). 
Accordingly $\eps_i=(y_i-x_i'\beta)/\sigma$ 
has mean value $e_i=(y_i-x_i'\hatt\beta)/\sigma$, 
the estimated residual, and variance $h_i^2={1\over n}x_i'M^{-1}x_i$, 
giving 
$$\eqalign{
F_{n,B}(t)&=\E\,\{F(t)\midd \data\} \cr
&=\E\,\Bigl[{a\over a+n}\Phi(t)+{n\over a+n} 
	{1\over n}\sum_{i=1}^n I\{ (y_i-x_i'\beta)/\sigma \le t\} 
		\midd \data\Bigr] \cr 
&={a\over a+n}\Phi(t)+{n\over a+n} 
	{1\over n}\sum_{i=1}^n\Phi\Bigl(
		{t-e_i\over h_i}\Bigr). \cr}$$
In contrast to the i.i.d.~case, see (1.5), 
this is a continuous distribution with density
$$f_{n,B}(t)={a\over a+n}\phi(t)+{n\over a+n}f_n(t)
	={a\over a+n}\phi(t)+{n\over a+n}
	{1\over n}\sum_{i=1}^n
	\phi\Bigl({t-e_i \over h_i}\Bigr){1\over h_i}. \eqno(7.4)$$
The second term $f_n(t)$ is a variable kernel density estimate
with smoothing parameters $h_i$ smaller than the usual ones, 
i.e.~$f_n(t)$ follows the ups and downs of a fine histogram more than
a typical kernel estimate would do. 
This particular result is implicit in Hjort (1987) and has also
been found by Olaf Bunke (1988). 

We can now describe the BB strategy. For a general prior $\pi(\beta,\sigma)$,
work out the posterior $\pi_n(\beta,\sigma\midd \data)$ and the 
corresponding generalisation of (7.4), with a $f_n(t)$ 
that is potentially more complicated but still a variable kernel 
estimate for estimated residuals 
$e_i=(y_i-x_i'\tilda\beta)/\tilda\sigma$. 
Choose a random $(\beta^*,\sigma^*)$ from $\pi_n(.)$ and then a BB sample 
$\eps^*_1,\ldots,\eps^*_{n+a}$ of size $n+a$ from 
$F_{n,B}(t\midd \beta^*,\sigma^*)=\E\,\{F(t)\midd \beta^*,\sigma^*,\data\}$, 
cf.~(7.3). Then compute the BB value 
$\theta^*_{BB}=\theta(\beta^*,\sigma^*,F^*_{BB})$, where
$F^*_{BB}$ is the empirical distribution of the chosen $\eps^*_i$'s. 
This is repeated a large number of times and gives
$\hatt G_n(t)={\rm Pr}_*\{\theta^*_{BB}\le t\}$, proposed 
as an approximation to 
$G_n(t)={\rm Pr}_B\{\theta(\beta,\sigma,F)\le t\midd \data\}$.  
BB-based point estimates and BB percentile intervals can then be computed.
Bias corrections of some sort can be carried out using the exact
information in (the parallel to) (7.4). 

An interesting bootstrapping strategy emerges in the 
case of vague prior information. This would mean 
a flat prior for $\beta$, a flat prior for $\log\sigma$, 
and $a\rightarrow0$ for the Dirichlet. The steps above 
take this form: Draw first $\sigma^*$ from the distribution that
corresponds to $1/\sigma^2$ being Gamma with parameters 
$\half(n-p)$ and $\half(n-p)\hatt\sigma^2$, where 
$\hatt\sigma^2=\sum_{i=1}^n(y_i-x_i'\hatt\beta)^2/(n-p)$ is 
the usual estimate. Then draw $\beta^*$ from 
$N\{\hatt\beta,{1\over n}(\sigma^*)^2M^{-1}\}$, and then 
$\eps^*_1,\ldots,\eps^*_n$ from the empirical distribution of 
$\eps_i=(y_i-x_i'\beta^*)/\sigma^*$. Finally compute
$\theta^*=\theta(\beta^*,\sigma^*,\eps^*_1,\ldots,\eps^*_n)$ 
based on these, i.e., based on pairs $(x_i,y_i^*)$
where $y_i^*=x_i'\beta^*+\sigma^*\eps^*_i$. 
This constitutes an alternative frequentist 
way of bootstrapping in the semiparametric regression model. 


\bigskip
{\bf 8. Bootstrapping schemes in hazard rate models.} 
To what extent do methods and results of the previous sections
generalise to situations with censored data, 
and to more general models for survival data analysis? 
This section deviates somewhat from the rest of the article 
and reports on a brief investigation into 
frequentist and Bayesian bootstrapping schemes 
for such problems, where it is natural to shift attention 
from cumulative distribution functions (c.d.f.'s)
to cumulative hazard rates (c.h.r.'s). 

\smallskip
{\sl 8A. From c.d.f.~$F$ to c.h.r.~$A$.} 
For concreteness we concentrate on the random censorship model
here. Generalisations to counting process models should not be difficult. 
Life-times $X_i^0$ from a distribution $F$ on $[0,\infty)$ may
be censored on the right, so that only $X_i=\min\{X_i^0,c_i\}$
and $\delta_i=I\{X_i\le c_i\}$ are observed. 
It is assumed that $X_i^0$'s and censoring times $c_i$'s are independent,
and that the $c_i$'s come from some $H$. 
The c.h.r.~$A$ is defined via $A[s,s+ds]=F[s,s+ds]/F[s,\infty)$, 
or $\dd A(s)=\dd F(s)/F[s,\infty)$ for short, which leads to 
$$A(t)=\int_0^t{\dd F(s)\over F[s,\infty)}
	\quad {\rm and} \quad 
	F(t)=1-\prod_{[0,t]}\{1-\dd A(s)\}. \eqno(8.1)$$
When $F$ is continuous then $A=-\log(1-F)$, 
but we will encounter non-continuous c.d.f.'s and c.h.r.'s,
for which the product integral representation (8.1) is appropriate;
see for example Hjort (1990). 
Parameters defined in terms of $F$ can equally be represented as 
functions of $A$, say $\theta=\theta_{\rm cdf}(F)=\theta_{\rm chr}(A)$. 

Introduce $N_n(t)=\sum_{i=1}^nI\{X_i\le t,\delta_i=1\}$, 
counting the number of observed events in $[0,t]$, 
and $Y_n(t)=\sum_{i=1}^nI\{X_i\ge t\}$,
the number at risk just before time $t$. The Kaplan--Meier 
and the Nelson--Aalen estimator are respectively 
$$F_n(t)=1-\prod_{[0,t]}\{1-\dd N_n(s)/Y_n(s)\}
	\quad {\rm and} \quad 
	A_n(t)=\int_0^t\dd N_n(s)/Y_n(s). \eqno(8.2)$$
Here $\dd N_n(s)$ jumps only at observed life-times, 
with jump $\Delta N_n(x_i)=1$ if these are distinct. 
In particular $A_n$ is the c.h.r.~associated with $F_n$, 
and $\dd A_n(s)=\dd N_n(s)/Y_n(s)$. In the uncensored case 
$\Delta A_n(x_i)=1/(n-i+1)$, assuming $x_1<\ldots<x_n$,
and then $F_n$ becomes the usual empirical c.d.f. 
Traditional nonparametric inference is based on the fact that
for large $n$, $A_n(.)$ has approximately independent increments with  
$$\E\,\dd A_n(s)\doteq\dd A(s), \quad 
{\rm Var}\,\dd A_n(s)=\E\,Y_n(s)^{-1}\,\dd A(s)\{1-\dd A(s)\}. \eqno(8.3)$$
A precise large-sample statement is that 
${\cal L}[\sqrt{n}\{A_n(.)-A(.)\}]\rightarrow V(.)$, 
in which $V(.)$ is a Gau\ss ian martingale with independent
increments and ${\rm Var}\,\dd V(s)=\dd A(s)\{1-\dd A(s)\}/y(s)$. 
Here $y(s)$ is the limit in probability of $Y_n(s)/n$, 
that is $y(s)=F[s,\infty)G[s,\infty)$ under present circumstances. 
See Hjort (1991), for example. 
In the continuous case $\dd A(1-\dd A)=\dd A$, of course. 

\smallskip
{\sl 8B. The weird bootstrap.} 
Let us make the following introductory remark,
which applies to both frequentist and Bayesian bootstrapping: 
There is nothing particularly magic about resampling data
per se, and other data-conditional simulation schemes 
might easily work as well. 
In the classical i.i.d.~framework, for example, 
resampling from $F_n$ creates a $F_n^*$ with the properties 
$$\E_*\,F_n^*(t)=F_n(t)
	\quad {\rm and} \quad 
	{\rm cov}_*\{F_n^*(s),F_n^*(t)\}=n^{-1}F_n(s)\{1-F_n(t)\}
	{\rm\ for\ }s\le t, \eqno(8.4)$$
which properly match 
$$\E\,F_n(t)=F(t)
	\quad {\rm and} \quad 
	{\rm cov}\{F_n(s),F_n(t)\}=n^{-1}F(s)\{1-F(t)\}
		{\rm\ for\ }s\le t. $$
This almost suffices for ${\cal L}_*\{\sqrt{n}(F_n^*-F_n)\midd \data\}$
to be close to ${\cal L}\{\sqrt{n}(F_n-F)\}$, and (3.7)--(3.8) make
this precise. 
But other simulation schemes that in one way or another 
create some artificial $F_n^*(.)$ with properties like (8.4) 
can also be expected to work. That is, even if the random 
$F_n^*$ is created from other means than actual sampling, 
one would expect 
${\cal L}_*[\sqrt{n}\{\theta_{\rm cdf}(F_n^*)
	-\theta_{\rm cdf}(F_n)\}\midd \data]$ 
to be close to 
${\cal L}[\sqrt{n}\{\theta_{\rm cdf}(F_n)-\theta_{\rm cdf}(F)\}]$.

In view of this remark and of (8.3) 
we should look for data-conditional simulation 
strategies that produce some random c.h.r.~$A_n^*$ 
with approximately independent increments and with 
$$\E_*\,\{\dd A_n^*(s)\midd \data\}\doteq \dd A_n(s), \quad 
{\rm Var}_*\{\dd A_n^*(s)\midd \data\}
	\doteq Y_n(s)^{-1}\,\dd A_n(s)\{1-\dd A_n(s)\}. \eqno(8.5)$$
Such schemes will succeed in the required 
$${\cal L}_*\bigl[\sqrt{n}\{\theta_{\rm chr}(A_n^*)
	-\theta_{\rm chr}(A_n)\}\midd \data\bigr]
	\doteq
  {\cal L}\bigl[\sqrt{n}\{\theta_{\rm chr}(A_n)
	-\theta_{\rm chr}(A)\}\bigr], \eqno(8.6)$$ 
or ${\cal L}_*\{\sqrt{n}(\hatt\theta_n^*-\hatt\theta_n)\midd \data\}
\doteq {\cal L}\{\sqrt{n}(\hatt\theta_n-\theta)\}$ for short, 
with few extra requirements. 
One very simple way of achieving this is to 
let $A_n^*(.)$ have independent increments and 
$$\dd A_n^*(s)=Y_n(s)^{-1}\,{\rm Bin}\{Y_n(s),\,\dd A_n(s)\}. \eqno(8.7)$$
So $A_n^*(.)$ is flat between observed life times, and 
at such a point $x_i$, say, the hazard jump  
$\Delta A_n^*(x_i)$ is a relative frequency from a 
binomial with parameters $Y_n(x_i)$ and $1/Y_n(x_i)$. 
This is Richard Gill's `weird bootstrap' (1990, personal communication). 

Note that $A_n^*$ corresponds to a random 
$F_n^*(x_i)=1-\prod_{x_j\le x_i}\{1-\Delta A_n^*(x_j)\}$,
which is different from that obtained through resampling from $F_n$. 
The weird bootstrap does not resample any data set, but it works,
with and without censoring. 
A precise asymptotic result about (8.6) can be proved.
In particular the weird bootstrap can be seen as an
alternative to Efron's bootstrap in the uncensored case,
developed from the hazard rate point of view. 

\smallskip
{\sl 8C. Exact nonparametric Bayesian analysis.}
Now we can embark on Bayesian issues. 
The canonical analogue to a Dirichlet process for $F$ is a
Beta process for $A$. Let $A$ be such a process with parameters 
$c(.)$ and $A_0(.)$, which means that $A$ has independent increments
that are approximately Beta distributed, 
$$\dd A(s)\approx {\rm Beta}\bigl[c(s)\dd A_0(s),\,
  c(s)\{1-\dd A_0(s)\}\bigr]. $$
Note that 
$$\E_B\,\dd A(s)=\dd A_0(s) \quad {\rm and} \quad 
	{\rm Var}_B\,\dd A(s)={\dd A_0(s)\{1-\dd A_0(s)\}\over c(s)+1}, $$
so $A_0$ is the prior guess and 
$c(s)$ is related to the concentration of the prior measure 
around this prior guess. When $c(s)=aF_0[s,\infty)$,
where $F_0=1-\prod_{[0,\cdot]}(1-\dd A_0)$,  
then $F=1-\prod_{[0,\cdot]}(1-\dd A)$ is Dirichlet with parameter $aF_0$. 
See Hjort (1990) about further properties for Beta processes.   

Given data $A$ is still a Beta process, with 
parameters $c+Y_n$ and $A_{n,B}$, where 
$$A_{n,B}(t)=\int_0^t{c(s)\,\dd A_0(s)+\dd N_n(s)\over c(s)+Y_n(s)} $$
is the Bayes estimate. 
So $A$ given data has independent increments with 
$$\dd A(s)\midd \data\approx {\rm Beta}
	\bigl[c(s)\,\dd A_0(s)+\dd N_n(s),
		\,c(s)\{1-\dd A_0(s)\}+Y_n(s)-\dd N_n(s)\bigr]. \eqno(8.8)$$
In particular 
$$\,E_B\,\{\dd A(s)\midd \data\}=\dd A_{n,B}(s)
	\quad {\rm and} \quad 
  {\rm Var}_B\{\dd A(s)\midd \data\}=
		{\dd A_{n,B}(s)\{1-\dd A_{n,B}(s)\}
		\over c(s)+Y_n(s)+1}. \eqno(8.9)$$
Full Bayesian posterior analysis is now theoretically possible, 
via simulation
of the independent increment process $A$ and then calculation of 
$\theta_{\rm chr}(A)$, leading in the end to 
$G_n(t)={\rm Pr}_B\{\theta_{\rm chr}(A)\le t\midd \data\}$, cf.~(1.3).
And in view of (8.3) and (8.9) we would get the pleasing result
$${\cal L}_B\bigl\{\sqrt{n}(A-A_{n,B})\midd \data\bigr\}
	\approx 
  {\cal L}\bigl\{\sqrt{n}(A_n-A)\bigr\},$$
i.e.~frequentists and Bayesians would agree for large sample sizes.
The full Bayesian simulation method is cumbersome, 
however, and requires a fine partitioning 
of the {half}line. 

\smallskip
{\sl 8D. A weird Bayesian bootstrap.} 
In view of the relative complexity of the full 
Bayesian Beta process approach one may look for 
approximating BB-strategies,
perhaps generalising our basic BB of (1.5)--(1.7). 
One way is to sample $X_1^{0*},\ldots,X_{n+a}^{0*}$ from the
Bayes estimate
$$F_{n,B}(t)=1-\prod_{[0,t]}\Bigl\{1-{c\,\dd A_0+\dd N_n\over c+Y_n}\Bigr\},$$
then pairing them with simulated censoring times $c_i^*$ from the
Kaplan--Meier estimate $H_n$ for $H$, and then
treating $X_i^*=\min\{X_i^{0*},c_i^*\}$ and 
$\delta_i^*=I\{X_i^{0*}\le c_i^*\}$ as a new BB-data set.
Here $a$ is related to the $c(.)$ function, for example being taken to
be its maximum value. 
This works, asymptotically, under some conditions, 
but not particularly well under non-negligible censoring. 
So in this sense there does not seem to be a natural generalisation
of this article's full data BB 
to hazard rate models with Beta process priors. 
The approximation suggested here does however work best
when $c(s)=aF_0[s,\infty)$, which is the ${\rm Dirichlet}(aF_0)$ 
prior for $F$, and indeed with BB sample size $n+a$. 

A simpler second solution which both works better 
and has a nice interpretation of its own 
is to generate $A^*_{BB}$ with independent 
binomial frequencies increments 
$$\dd A^*_{BB}(s)=\{c(s)+Y_n(s)\}^{-1}\,{\rm Bin}\{c(s)+Y_n(s),
		\dd A_{n,B}(s)\}. \eqno(8.10)$$
This manages to almost match (8.9), and the small difference 
disappears asymptotically. 
At observed life times the jump $\Delta A^*_{BB}(x_i)$ 
is a binomial $[c(x_i)+Y_n(x_i),1/\{c(x_i)+Y_n(x_i)\}]$
divided by $c(x_i)+Y_n(x_i)$. Again, this scheme does not 
correspond to data resampling, but weirdly kills and reincarnates 
individuals at each time point. 

\smallskip
{\sl 8E. Exact Bayesian and BB analysis under a noninformative 
reference prior.}
Let $c(.)$ go to zero in the above constructions. 
The exact Bayesian solution is then a Beta process
$A$ with parameters $Y_n$ and $A_n$, i.e.
$$\dd A(s)\midd \data\sim {\rm Beta}\{\dd N_n(s),
   Y_n(s)-\dd N_n(s)\}, \eqno(8.11)$$
with independent jumps only at observed life times. 
Note that 
$$\E_B\,\{\dd A(s)\midd \data\}=\dd A_n(s) 
	\quad {\rm and} \quad 
  {\rm Var}_B\{\dd A(s)\midd \data\}
	={\dd A_n(s)\{1-\dd A_n(s)\} \over Y_n(s)+1}.$$
Observe also that letting $c(.)\rightarrow0$ in the posterior
distributions is the same as letting 
$a\rightarrow0$ in the posterior distribution with
a ${\rm Dirichlet}(aF_0)$. 
In this way we have arrived at a generalisation to censored data
for Rubin's noninformative Bayesian bootstrap. 
In addition to being a natural limit of proper Bayes 
solutions it can be given an empirical Bayesian interpretation. 
The (8.11) method is also the method proposed by Lo (1991);
see his paper for further properties. 

Letting $c(.)\rightarrow0$ in the weird BB of 8D 
gives Gill's weird bootstrap of 8B. 
The latter can therefore be seen as the noninformative limit
of a natural simulation-based approximation to a full Bayesian
method. 

One can prove that all the schemes described here are first order 
equivalent. In particular each scheme will reach confidence intervals 
asymptotically equivalent to (3.11).  

\smallskip
{\sl 8F. Cox regression.}
Let us finally note that the methods above can be extended 
and used in the semiparametric Cox regression model. 
Suppose individual no.~$i$ has covariate $z_i$ and 
c.h.r.~$A_i$, and that 
$$1-\dd A_i(s)=\{1-\dd A(s)\}^{\exp(\beta z_i)}, 
	\quad i=1,\ldots,n.$$ 
If the Bayesian prior is that $\beta$ comes from some $\pi(\beta)$ 
and that $A$ independently is a Beta process $(c,A_0)$,
then the posterior distributions can be worked out,
making a full semiparametric Bayesian analysis awkward but possible,
through cumbersome simulations. 
See Hjort (1990, Section 6).
Simulation-based approximations to this scheme can be developed,
with ideas as above, giving in particular a weird BB scheme,
but requiring more involved distributions than the simple binomial. 
Let us merely report on the noninformative limit version of 
the exact Bayes solution, as $c(.)\rightarrow0$. 
First draw a $\beta$ from 
$$\pi_n(\beta)={\rm const.}\,
	\prod_{i=1}^n \bigl\{\psi\bigl(R_n(x_i,\beta)\bigr)
     	-\psi\bigl(R_n(x_i,\beta)-\exp(\beta z_i)\bigr)\bigr\}^{\delta_i}.$$
Then let $A$ be flat between observed life times, 
and have independent jumps 
$$\Delta A(x_i)\sim 
 	{z^{-1}\bigl\{(1-z)^{R_n(x_i,\beta)-\exp(\beta z_i)-1}
		-(1-z)^{R_n(x_i,\beta)-1}\bigr\}
	\over 
	\psi\bigl(R_n(x_i,\beta)\bigr)
		-\psi\bigl(R_n(x_i,\beta)-\exp(\beta z_i)\bigr)},
		\quad 0<z<1, $$
for those $x_i$ with $\delta_i=1$. 
In these expressions $R_n(s,\beta)=\sum_{i=1}^n Y_i(s)\exp(\beta z_i)$,
where $Y_i(s)=I\{X_i^0\ge s,c_i\ge s\}$ is the at-risk indicator 
for individual $i$. And any sensible simpler way of simulating
$\Delta A^*_{BB}(x_i)$ instead, with asymptotically correct
matching for the two first moments, defines a weird BB. 
		
\bigskip
{\bf 9. Supplementing results and remarks.} 
This final section offers some concluding comments and mentions
some extensions of previous results. 

\smallskip
{\sl 9A. Two viewpoints.} 
There are presumably two ways to approach the problem 
of handling $\theta(F)$ in a Bayesian nonparametric way.
One way is to ignore the underlying $F$ and 
concentrate on $\theta(F)$ and what the prior information on 
this particular parameter is. In the end some Bayesian
calculations are carried out for $\theta$ given data.
In this mode each parameter must be treated separately, 
and inconsistencies can occur, since Bayesians are nonperfect. 
The second way is the one chosen in this article,
where information is expressed in terms of the underlying $F$
once and for all, after which analysis can proceed on an
automatic basis for every conceivable $\theta(F)$. 

\smallskip
{\sl 9B. Exact simulation.}
There are actually ways of simulating almost exactly from 
$G_n={\cal L}_B\{\theta(F)\midd \data\}$,
cf.~remarks made at the start of 2A,
where one such method was described, using a fine partition of the
real line and finite-dimensional Dirichlet distributions. 
Another way would be through simulation of $F$ via its
product integral representation in terms of the cumulative 
hazard process $A$, which is a Beta process, 
see Section 8 above. This remark and results there 
show that such posterior simulation of $\theta(F)$ is possible
even with censored data, and in more complex models 
for survival data. 

A third possibility is to use Sethuraman's 
constructive definition of an arbitrary Dirichlet process, 
see Sethuraman and Tiwari (1982). 
The present ${\rm Dir}(aF_0+nF_n)$ case (see (1.2)) 
can be represented as follows. 
Generate an infinite i.i.d.~sequence $\{x_i'\}$ from $F_{n,B}$ of (1.5),
and an infinite i.i.d.~sequence $\{B_i\}$ from ${\rm Beta}\{1,a+n\}$.
Then let $A_i=B_i\prod_{j=1}^{i-1}(1-B_j)$ and use 
$$F=\sum_{i=1}^\infty A_i\delta(x_i'), \eqno(9.1)$$
where $\delta(x)$ means unit point mass at position $x$.  
To see how this can be used, 
consider the mad-parameter $\theta(F)=\int|x-{\rm med}(F)|\,\dd F(x)$,
for example. Approximate $F$ by using a large number $I$ instead of 
$\infty$ in (9.1), perhaps $I=1000$. Order the $x_i'$ points and 
determine the one for which the cumulative probability mass
$\sum_{i=1}^jA_i$ first exceeds $\half$; this gives an approximation
${\rm med}'$ to $F^{-1}(\half)$. 
Some care is required since there will be heavy ties in the $x_i'$ data. 
Then compute $\theta'=\sum_{i=1}^I|x_i'-{\rm med}'|A_i$, all in all
giving an approximation to one particular $\theta$ drawn from $F$.
This algorithm must then be repeated a large number of times to 
form ${\cal L}\{\theta(F)\midd \data\}$. --- This elaborate
strategy makes almost-exact Bayesian caculations possible, 
and in a certain sense makes the BB less necessary. 
But arguments still favouring the BB include 
(i) that it is much simpler
to use, regarding both programming, simulation, and cpu-use, 
(ii) that it is less functional-dependent,
(iii) that BB and almost-exact simulation are first 
order equivalent, by Section 3, and 
(iv) that the BB perhaps is more trustworthy and realistic 
than the almost-exact version in that it only exploits the 
first and second order characteristics of the Dirichlet process
structure, and not the more esoteric ones, like the inherent 
discreteness of its sample paths, visible in (9.1). 
In any case the (9.1)-based method does make 
almost-exact posterior Dirichlet analysis possible
and should be included in any serious comparison between 
the various strategies. 

\smallskip
{\sl 9C. Invariance under transformations.}
Both the nonparametric Bayesian confidence interval 
$[\theta_L,\theta_U]$ of (1.4) and its BB approximation
$[\hatt\theta_L,\hatt\theta_U]$ of (1.9) transform very neatly,  
with respect to both data-transformations and parameter-transformations.
(i) Suppose $\nu=g(\theta)$ is a new parameter, with a smooth
and increasing $g(.)$. The (1.4) scheme uses 
$H_n(t)={\rm Pr}_B\{\nu(\theta(F))\le t\}=G_n(g^{-1}(t))$,
and the (1.9) uses $\hatt H_n(t)={\rm Pr}_*\{\nu^*_{BB}\le t\}
	=\hatt G_n(g^{-1}(t))$. It follows that 
$$[\nu_L,\nu_R]=[g(\theta_L),g(\theta_R)],
	\quad 
  [\hatt\nu_L,\hatt\nu_R]=[g(\hatt\theta_L),g(\hatt\theta_R)]. \eqno(9.2)$$
(ii) Suppose $Y_i=h(X_i)$ for a smooth increasing $h(.)$. 
If $F$ for $X_i$ is Dirichlet $aF_0$, then $\tilda F=Fh^{-1}$ 
for $Y_i$ is Dirichlet $a\tilda F_0=aF_0h^{-1}$. 
Write $\nu(\tilda F)=\theta(F)$ for the old parameter seen 
in the context of $Y_i$'s from $\tilda F$. 
Then $H_n(t)={\rm Pr}_B\{\nu\le t\}=G_n(t)$ and 
the Dirichlet transformation property implies 
$\hatt H_n(t)={\rm Pr}_*\{\nu^*_{BB}\le t\}=\hatt G_n(t)$. 
So (1.4) and (1.9) are invariant under data transformations. 

\smallskip
{\sl 9D. Boot sample size.}
The bootstrap sample size `${\rm boot}$' in (1.8) should of course
be large in order for (1.8) and (1.9) to work well,
i.e.~for functions of $\hatt G_{n,{\rm boot}}(.)$ to be close to
the same functions of $\hatt G_n$. 
The investigation of Efron (1987, Section 9),
albeit for a different bootstrap, 
is relevant here, and indicates that 
${\rm boot}=1000$ may be a rough minimum for quantiles in the tail,
required in (1.9), but that ${\rm boot}=100$ may suffice 
for average operations like the mean, required in (1.8).  

\smallskip
{\sl 9E. Highest posterior density.} 
The starting point for our confidence intervals has been (1.4).
Sometimes in the Bayesian literature 
highest posterior density regions are advocated instead.
In the present case this would involve approximating the posterior
distribution $G_n$ with one with a density $g_n$, and then
letting $\{t\colon g_n(t)\ge g_0\}$ be the confidence region,
for appropriate level $g_0$. This approach makes most sense when 
$g_n$ is unimodal, which it would not necessarily be in applications
of the present kind, due to the fact that the posterior distribution
of $F$ places extra weight on the observed data points. 
This is illustrated in Section 5 for the case of the median.

\smallskip
{\sl 9F. Data-dependent functionals.} 
The functional $\theta=\theta(F)$ can depend on the sample size;
the described BB procedure works specifically for the given $n$.
It is also allowed to depend upon the 
actual data sample, say $\theta=\theta(F,x_1,\ldots,x_n)$.
Let us illustrate this comment with a description of how a nonparametric
Bayesian might construct a simultaneous confidence band for $F$. 
Consider 
$$\theta_{\min}=\min_{a\le t\le b}{F(t)-F_{n,B}(t)\over 
	[F_{n,B}(t)\{1-F_{n,B}(t)\}]^{1/2}}
	\quad {\rm and} \quad 
  \theta_{\max}=\max_{a\le t\le b}{F(t)-F_{n,B}(t)\over 
	[F_{n,B}(t)\{1-F_{n,B}(t)\}]^{1/2}}.$$
The natural band is
$$F_{n,B}(t)-c[F_{n,B}(t)\{1-F_{n,B}(t)\}]^{1/2}
	\le F(t) \le
  F_{n,B}(t)+d[F_{n,B}(t)\{1-F_{n,B}(t)\}]^{1/2} $$
for $a\le t\le b$, where $c$ and $d$ 
ideally would be determined by 
$${\rm Pr}_B\{-c\le \theta_{\min}(F,x_1,\ldots,x_n),\,
	        \theta_{\max}(F,x_1,\ldots,x_n)\le d\,\midd \data\}=0.90,$$
say (with an additional condition to make them unique, like
requiring minimisation of $c+d$). 
The BB method consists of generating perhaps 1000 values of
$$\theta^*_{\min,BB}=\min_{a\le t\le b}{F_{BB}^*(t)-F_{n,B}(t)\over 
	[F_{n,B}(t)\{1-F_{n,B}(t)\}]^{1/2}}, 
	\,\,
  \theta^*_{\max,BB}=\max_{a\le t\le b}{F_{BB}^*(t)-F_{n,B}(t)\over 
	[F_{n,B}(t)\{1-F_{n,B}(t)\}]^{1/2}},$$
and using the correspondingly defined $\hatt c$ and $\hatt d$. 
One may prove that $(n+a)^{1/2}(\hatt c-c)$ 
and $(n+a)^{1/2}(\hatt d-d)$ go to zero in probability,
by methods and results of Section 3. 
[Strictly speaking, this is true provided 
${\rm boot}_n$ realisations are generated instead of 1000 
and ${\rm boot}_n/(n\log n)$ grows with towards infinity.]
It could be advantageous to use this asymmetric band 
instead of the simpler symmetric one 
since the distribution of $F-F_{n,B}$ is typically skewed. 

\smallskip
{\sl 9G. Two-sample BB.} 
The BB method can be generalised to two-sample situations,
and indeed to more general non-i.i.d.~models, as shown in 
Section 7. To illustrate,
let $x_1,\ldots,x_n$ and $y_1,\ldots,y_m$ be samples from 
respectively $F_1$ and $F_2$, assume 
$\theta=F_1^{-1}(\half)-F_2^{-1}(\half)$ is of interest, 
and suppose $F_1\sim {\rm Dir}(aF_{1,0})$ and 
$F_2\sim {\rm Dir}(bF_{2,0})$. 
A Bayes estimate and confidence interval for this difference of
population medians can be obtained by generating perhaps 1000
realisations of 
$\theta_{BB}^*={\rm med}\{X^*_1,\ldots,X^*_{n+a}\}
	-{\rm med}\{Y^*_1,\ldots,Y^*_{m+b}\}$,
where the $X^*_i$'s are drawn from $(aF_{1,0}+nF_{1,n})/(a+n)$
and the $Y^*_i$'s from $(bF_{1,0}+mF_{2,m})/(b+m)$,
and then treating the resulting histogram (or smoothed density
estimate) as the posterior distribution of $\theta$. 

\bigskip
\centerline{\bf References}

\smallskip
\parindent20pt
\parskip2pt

\def\ref#1{%
  \par
  \hangindent=\parindent
  \hangafter=1
  \noindent #1\par
}


\ref{%
Bickel, P.J.~and Freedman, D.A. (1981).
Some asymptotic theory for the bootstrap.
{\sl Annals of Statistics}~{\bf 9}. 1196--1217.} 

\ref{%
Billingsley, P. (1968).
{\sl Convergence of Probability Measures.}
Wiley, Singapore.}

\ref{%
Boos, D.D.~and Serfling, R.J. (1980).
A note on differentials and the CLT and LIL for 
statistical functions with application to M-estimates.
{\sl Annals of Statistics}~{\bf 8}, 618--624. }

\ref{%
Bunke, O. (1988).
Posterior distributions in semiparametric models.
Technical report, Humboldt Universit\"at Berlin.} 

\ref{%
Cifarelli, D.M.~and Regazzini, E. (1990).
Distribution functions of means of a Dirichlet process.
{\sl Annals of Statistics}~{\bf 18}, 429--442.} 

\ref{%
Efron, B. (1979). 
Bootstrap methods: Another look at the jackknife.
{\sl Annals of Statistics}~{\bf 7}, 1--26.}

\ref{%
Efron, B. (1982). 
{\sl The Jackknife, the Bootstrap, and Other Resampling Plans.}
CBMS 38, SIAM--NSF.}

\ref{%
Efron, B. (1987).
Better bootstrap confidence intervals (with discussion).
{\sl Journal of the American Statistical Association}~{\bf 82}, 171--200. }

\ref{%
Ferguson, T.S. (1973).
A Bayesian analysis of some nonparametric problems.
{\sl Annals of Statistics} {\bf 1}, 209--230.}

\ref{%
Ferguson, T.S. (1974).
Prior distributions on spaces of probability measures.
{\sl Annals of Statistics} {\bf 2}, 615--629.}

\ref{%
Gill, R.D. (1989).
Non- and semiparametric maximum likelihood estimation and
the von Mises method (part I, with discussion).
{\sl Scandinavian Journal of Statistics}~{\bf 16}, 97--128.} 

\ref{%
Gill, R.D. (1990).
A weird bootstrap.
[Unpublished; personal communication.]}

\ref{%
Hannum, R., Hollander, M., and Langberg, N. (1981).
Distributional results for random functionals of
a Dirichlet process.
{\sl Annals of Probability}~{\bf 9}, 665--670.  }

\ref{%
Hjort, N.L. (1976).
Applications of the Dirichlet process to some nonparametric 
problems (in Norwegian). Graduate thesis, University of Troms\o.
[Abstract in {\sl Scandinavian Journal of Statistics}~{\bf 4}, 1977, p.~94.]}

\ref{%
Hjort, N.L. (1985).
Bayesian nonparametric bootstrap confidence intervals.
NSF- and LCS-Technical report,
Department of Statistics, Stanford University. }

\ref{%
Hjort, N.L. (1986).
Contribution to the discussion of Diaconis and Freedman's
`On the consistency of Bayes estimates'.
{\sl Annals of Statistics}~{\bf 14}, 49--55. }

\ref{%
Hjort, N.L. (1987).
Semiparametric Bayes estimators.
Proceedings of the First World Congress of the Bernoulli Society,
Tashkent, USSR, 31--34. VNU Science Press.}  

\ref{%
Hjort, N.L. (1988). 
Contribution to the discussion of Hinkley's lectures on bootstrapping.
To appear in {\sl Scandinavian Journal of Statistics} }

\ref{%
Hjort, N.L. (1990).
Nonparametric Bayes estimators based on Beta processes in models
for life history data.
{\sl Annals of Statistics}~{\bf 18}, 1259--1294. }

\ref{%
Hjort, N.L. (1991).
Semiparametric estimation of parametric hazard rates.
Proceedings of the 
{\sl NATO Advanced Study Workshop on 
Survival Analysis and Related Topics}, 
Columbus, Ohio, June 1991. }

\ref{%
Lo, A.Y. (1987).
A large-sample study of the Bayesian bootstrap.
{\sl Annals of Statistics}~{\bf 15}, 360--375.}

\ref{%
Lo, A.Y. (1991). 
A Bayesian bootstrap for censored data.
Technical report, Department of Statistics, SUNY at Buffalo. }

\ref{%
Newton, M.A.~and Raftery, A.E. (1991). 
Approximate Bayesian inference by the weighted likelihood bootstrap.
Technical report, Department of Statistics, University of Washington, 
Seattle.}

\ref{%
Parr, W.C. (1985).
The bootstrap: some large sample theory, and connections with robustness.
{\sl Statistics and Probability Letters}~{\bf 3}, 97--100. }

\ref{%
Rubin, D.B. (1981). The Bayesian bootstrap.
{\sl Annals of Statistics}~{\bf 9}, 130--134. }

\ref{%
Sethuraman, J.~and Tiwari, R. (1982).
Convergence of Dirichlet measures and the interpretation of their
parameter. In {\sl Proceedings of the Third
Purdue Syposium on Statistical Decision Theory and Related Topics},
S.S.~Gupta and J.~Berger (eds.), 305--315. Academic Press, New York. }

\ref{%
Shao, J. (1989).
Functional calculus and asymptotic theory for statistical analysis.
{\sl Statistics and Probability Letters}~{\bf 8}, 397--405.}

\ref{%
Singh, K. (1981).
On the asymptotic accuracy of Efron's bootstrap.
{\sl Annals of Statistics}~{\bf 9}, 1187--1195. }

\ref{%
Yamato, H. (1984).
Characteristic functions of means of distributions
chosen from a Dirichlet process. 
{\sl Annals of Probability}~{\bf 12}, 262--267. }

\bigskip
\parindent2.0truecm
\baselineskip11pt
\parskip0pt
\obeylines
Nils Lid Hjort 
Department of Mathematics and Statistics
University of Oslo
P.B.~1053 Blindern 
N--0316 Oslo 3, Norway
e-mail: {\tt nils@math.uio.no}

\bye